\documentclass[absolute]{mymathart}
\usepackage[theorems]{mymathmacros}

\usepackage{subfigure,graphicx}
\usepackage{amsmath,amssymb}

\title{On Prime Ends and Local Connectivity}
\author{Lasse Rempe-Gillen}
\address{Department of Mathematical Sciences \\
    University of Liverpool \\
     Liverpool, United Kingdom L69 7ZL}
\email{l.rempe@liverpool.ac.uk}
\thanks{Supported in part by EPSRC grant EP/E017886/1.}
\dedicatory{Dedicated to the memory of Professor Gerald Schmieder}
\date{\today}

\subjclass[2000]{Primary 30C35; Secondary 54F15, 30D05, 37F10}

\begin{document} 

 \begin{abstract}
  Let $U\subset\Ch$ be a simply connected domain whose
   complement $K=\Ch\setminus U$ contains more than one point.
   We establish a characterization of the local connectivity of $K$ at
   a point $z_0\in\partial U$ in terms of the prime ends of $U$ whose 
   impressions contain $z_0$.
   Invoking a result of Ursell and Young \cite{ursellyoung},
   we obtain an alternative proof of a theorem of Torhorst, which states 
   that
   the impression of a prime end of $U$ contains at most
   two points at which $K$ is locally connected. 
 \end{abstract}

 \maketitle

 \section*{Historical comment}
  This article appeared in the Bulletin of the
   London Mathematical Society in 2008. In 2012, Donald Sarason kindly
   pointed out to me that Theorem \ref{thm:main} was 
   proved by Marie Torhorst in 1918 in her dissertation at Bonn University. 
   Her original thesis has been lost (a copy exists neither at Bonn nor 
    in her personal \emph{Nachla{\ss}}), but 
    she published the result in \cite{torhorst}. 
    It seems that this theorem has largely been forgotten: 
    no reference to Torhorst's paper is recorded on MathSciNet, while 
    the latest reference that I was able to find on Zentralblatt
    dates from 1930. Whyburn \cite{whyburn} references her article
    in his book \emph{Analytic Topology}, but does not mention the theory
    of prime ends (see below). 

  In the 1960s, Sarason 
   himself wrote a paper \cite{sarason}
   with the exact same title as this present note, 
   which also proves Torhorst's result using the more recent work of Ursell 
   and Young. 
   However, his manuscript was not accepted for publication at the time, 
   as he explains:

\begin{quotation}
I submitted the paper to the Michigan Math. J., then edited by George
  Piranian, the person who taught me about prime ends and much more about 
  complex analysis. (George is one of my mathematical heroes.) George 
  discussed the paper with Collingwood, one of his collaborators. Their 
  conclusion was that interest in prime ends at the time was at such a low ebb 
  that the paper was likely to be largely ignored.

I did publish an abstract of the paper in the 
  Notices of the A.M.S. (Vol 16 (1969), p. 701). At the time the Notices 
  published abstracts of talks given at society meetings, plus what I think 
  were called by-title abstracts, which any member of the society could use 
  to announce a result. If my memory is correct, I received as a result of 
  the abstract only one request for a copy of the paper.
\end{quotation}
 
 As far as I am aware, Theorem \ref{thm:characterization}, which is 
  a characterization of local connectivity at a given point, and of which
  Theorem \ref{thm:main} is a corollary, has not previously appeared
  elsewhere. 
  (Note, however, that the argument that proves the ``only if''
   direction is the same as the one that appears already in \cite{sarason}, 
    which also contains the ``if'' direction in the special case that every 
    prime end whose impression contains the point in question is of the first 
    kind.)

 Sarason's manuscript \cite{sarason}, 
  along with George Piranian's letter and the announcement in the Notices, 
  are available from my professional website. 
  The preprint version of my original article will
  follow below, after some additional remarks concerning
  Marie Torhorst's mathematical work.

\section*{Further comments on Torhorst's work}

  Torhorst's paper \cite{torhorst} contains a number of other interesting results,
     including the fact that the impression of each prime end can contain
     at most three points at which the boundary is ``accessible from all sides'' in the
     sense of Schoenflies (compare \cite[p.\ 111]{whyburn} for the definition).
     Moreover, she appears to state and prove for the first time 
     that 
     local connectivity of the boundary of a domain is equivalent to the 
     continuous extension of the Riemann map. This result has come to be
     extremely prominent, particularly in the area of one-dimensional holomorphic
     dynamics, but is usually attributed to Carath\'eodory (compare Corollary \ref{cor:caratheodory} below). 
     While Carath\'eodory's theory
     of prime ends is of course central to the proof, to my knowledge
     he never made a connection to local connectivity. Indeed, it seems that the latter notion
     was only developed by Hahn at a time at which Carath\'eodory's theory was already completed. 
     In view of this fact, it would seem more appropriate to refer to the result as the 
      \emph{Carath\'eodory-Torhorst Theorem}. 

   Today, Torhorst is remembered, if at all, for what is sometimes called the \emph{Torhorst Theorem} (see 
     \cite[Theorem 2.2 in Chapter VI]{whyburn}): 
     the boundary of each complementary domain of a locally connected continuum is itself
     locally connected. Prof.\ Walter Purkert has kindly pointed out that Hausdorff also studied
     this theorem \cite{hausdorfftorhorst}. However, this is only a corollary of Torhorst's main
     results. It is unfortunate
     that her contributions to the theory of prime ends appear to 
     have been forgotten, particularly considering 
     the relatively small number of female researchers in mathematics at that time. 

   Torhorst did not become a professional research mathematician~-- her autobiography articulates a lack of
     confidence in her own mathematical ability, and in any case, academic research at the time was
     an almost exclusively male undertaking. (In 1919 Emmy Noether became the first woman to be allowed to to habilitate in mathematics in
     Germany, with the support of Hilbert and Klein but only after years of considerable resistance from the academic establishment
   in G\"ottingen.) 
   A life-long communist, Torhorst was persecuted by the Nazi regime and later became a functionary
   in the SED, the governing socialist party of East Germany. In 1947, she was appointed Minister
   of Education in the state of Th\"uringen, which 
   makes her the first-ever female minister 
    (on state or national level) in the history of Germany.
    Marie Torhorst died in May 1989 
    at the age of 100. More information about her life can be found at 
    \cite{torhorstwebsite} and in the sources referenced there.    

 \section{Introduction}
  The theory of \emph{prime ends} was developed by
   Carath\'eodory \cite{caratheodory} in 1913.
   One of its central theorems states that the complement $K$ of 
   a simply connected domain $U=\Ch\setminus K$, $\#K>1$, 
   is locally connected if and only
   if every prime end has trivial impression, which in turn 
   is equivalent to any
   Riemann map $\phi:\D\to U$ having a continuous extension to
   $\partial\D$. (See Section \ref{sec:primeends} for a short
   introduction to the standard definitions and terminology of prime end
   theory.) 
  This note investigates the question whether there is a relationship
   between local connectivity of $K$ \emph{at a point $z_0\in\partial U$}
   --- i.e., the existence of arbitrarily small connected neighborhoods
    of $z_0$ in $K$,\footnote{Sometimes this property is instead referred to as \emph{connected im kleinen}.}
 compare Section \ref{sec:localconnectivity} --- 
   and the structure of the prime ends of $U$ whose impressions
   contain $z_0$.

  This question 
   seems very natural, 
   and may be of particular interest due to the 
   prominence that local connectivity of Julia sets and the Mandelbrot
   set at certain points has received in recent years 
   (see e.g.\ \cite{hubbardyoccoz,kiwifibers}). 
   However, 
   it does not appear to have received any treatment in the 
   literature so far. 

 A naive hope might be that $K$ is locally connected at $z_0$ if and only
  if every prime end impression which contains $z_0$ is trivial, but
  this is false, as 
  the well-known case of the ``double comb''
   shows (Figure \ref{fig:doublecomb}).
  However, study of this and similar examples suggests that a nontrivial
  impression
  should not contain ``too many'' points of local connectivity. 
  In this note, we
  demonstrate 
  that ``not too many'' can be made very precise.
  In fact, the example in Figure \ref{fig:doublecomb} 
  is already best possible. 

 \begin{thm}[Prime ends and local connectivity]
   \label{thm:main}
  Let $U\subset\Ch$ be a simply connected domain such that
   $K := \Ch\setminus U$ contains more than one point, and let $p$ be a
   prime end of $U$. Then the impression $I(p)$ contains at most two
   points at which $K$ is locally connected. 
 \end{thm} 
 \begin{remark}  
  We also show that, if furthermore 
   the prime end $p$ is \emph{symmetric}
   (see Section \ref{sec:ursellyoung}), 
   then $I(p)$ contains at most 
   \emph{one} point at which $K$ is locally connected.
 \end{remark}

 The proof of 
  Theorem \ref{thm:main} 
  uses a result (Theorem \ref{thm:ursellyoung})
  concerning the ``wings'' (aka the left and right cluster sets) 
  of a prime end
  that was proved by
  Ursell and Young \cite{ursellyoung} in 1951 and
  deserves to be far better known.

 We will deduce Theorem \ref{thm:main} from Theorem
  \ref{thm:ursellyoung} by
  developing a necessary and sufficient criterion for
  local connectivity at $z_0$ in terms of the prime ends of $U$. 
  To state this result, we introduce the following notion.
  \begin{defn}[Strong minimality] \label{defn:strongminimality}
   Let $p$ be a prime end, and let $z_0$ belong to the
    impression of $p$. We  say that
    \emph{$z_0$ is strongly minimal in $p$} if,
    for every sequence $w_j\in U$ converging to $p$ which does not
    accumulate on $z_0$, there is a curve $\Gamma:[0,\infty)\to U$ that
    converges to $p$ and passes
    through
    all $w_j$  but does not accumulate on $z_0$ (as $t\to\infty$). 
  \end{defn}
  This terminology is motivated by such a point being minimal 
   with respect to Ursell and Young's ordering by priority; 
   see Definition \ref{defn:priority}.

 \begin{thm}[Characterization of local connectivity]
    \label{thm:characterization}
  Let $z_0\in \partial U$. Then $K$ is locally connected at $z_0$
   if and only if $z_0$ is strongly minimal in every prime end whose
   impression contains $z_0$. 
 \end{thm}
  
 The proof of Theorem \ref{thm:characterization} is elementary,
  and the result
  might almost be considered a restatement of the definition
  of local connectivity. However, it does provide an
  interesting and quite satisfying answer to our initial
  question about the connection between prime ends and local connectivity;
  in particular it contains Carath\'eodory's characterization
  of local connectivity of $K$ (see Corollary \ref{cor:caratheodory}).
  Theorem \ref{thm:main} follows from Theorem \ref{thm:characterization}
  and the aforementioned
  result by Ursell and Young (compare Corollary \ref{cor:stronglyminimal}).

 \subsection*{Basic notation}
  We denote the complex plane by $\C$, the Riemann sphere by $\Ch$, and
   the unit disk by $\D$. We 
   write $\D_{\delta}(z)$ for the (Euclidean) disk of
   radius $\delta$ around $z$.

 \subsection*{Organization of the article}
  In Section \ref{sec:localconnectivity}, we define local connectivity
   at a point and discuss
   a number of variations of this definition. We also
   develop a simple characterization
   of local connectivity of $K$ at $z_0$.
   Section \ref{sec:primeends} provides a short review of the theory
   of prime ends and the proof of Theorem \ref{thm:characterization}.
   In Section \ref{sec:ursellyoung}, we discuss Theorem
   \ref{thm:ursellyoung}, by 
   Ursell and Young, 
   and deduce Theorem \ref{thm:main} from it. For
   completeness, we provide a proof of Theorem \ref{thm:ursellyoung}
   in the Appendix. 

 \subsection*{Acknowledgments} 
  I had many interesting and enjoyable discussions on this subject over
   the years,
   in particular with Chris Bishop, David Epstein, 
   Christian Pommerenke, Lex Oversteegen, Dierk Schlei\-cher and the late
   Gerald Schmieder. 
   I would like to thank Walter Bergweiler for a choice of seminar topic that
   not only introduced me to
   Pommerenke's excellent book \cite{pommerenke}, but also led me to discover
   Theorem \ref{thm:main} as an undergraduate at Kiel University in 1999. 
   Finally, I am grateful to Christian Pommerenke and
   Lex Oversteegen for encouraging me to 
   publish this note.

\begin{figure}%
  \subfigure[The double comb\label{fig:doublecomb}]%
                     {\resizebox{.425\textwidth}{!}{%
                               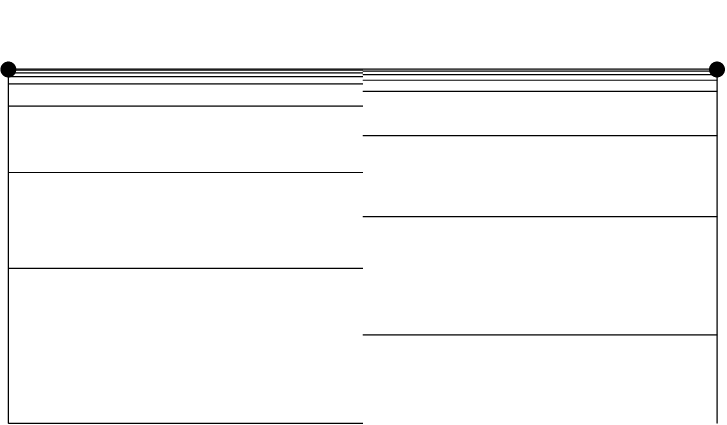}}%
    \hfill
  \subfigure[The witch's broom\label{fig:witchesbroom}]{%
                      \resizebox{.425\textwidth}{!}{%
                               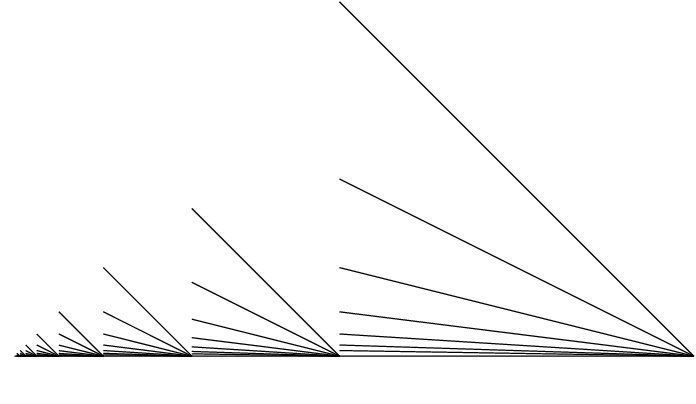}}\\%
  \subfigure[Path connectivity]{\label{fig:pathlc}%
                    \resizebox{.425\textwidth}{!}{%
                               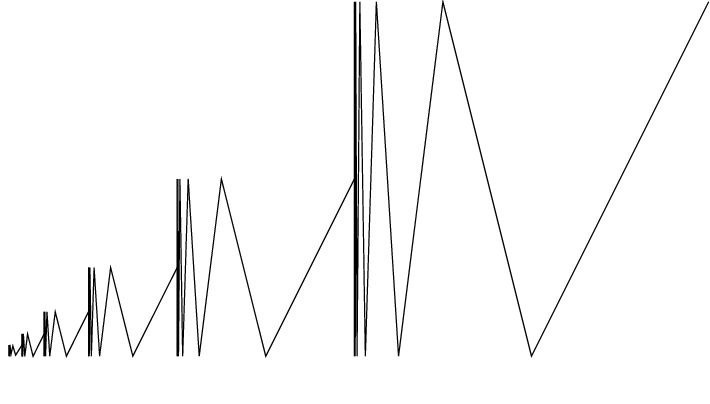}}%
    \hfill
  \subfigure[$K$ vs. $\partial U$\label{fig:boundary}]{%
                    \resizebox{.425\textwidth}{!}{%
                               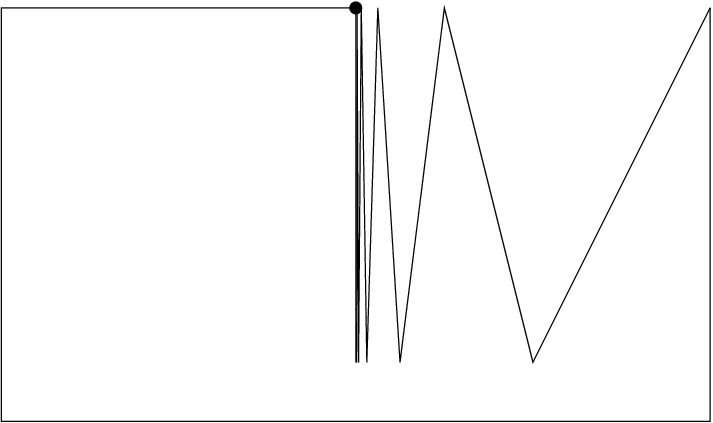}}\\%
 \caption{Several examples of simply connected domains and their boundaries.
   (a) illustrates Theorem \ref{thm:main}: the
     interval at the top of
     the figure is the impression of a single prime end $p$, and
     $K=\Ch\setminus U$ is locally connected at the two endpoints $z_0$
     and $z_1$. Also note that these endpoints are \emph{strongly minimal}
     in the sense of Definition \ref{defn:strongminimality}, while the interior
     points are not. 
 Examples (b) to
     (d) illustrate our remarks on the
     definition of local connectivity in Section \ref{sec:localconnectivity}.%
     \label{fig:examples}}
\end{figure}

\section{Local connectivity} \label{sec:localconnectivity}
  For the remainder of the paper, let $U\subset\Ch$ 
   be a simply connected domain
   whose complement
    $K := \Ch\setminus U$ contains at least two points. 
    We will assume without loss of generality
    that $\infty\in U$, so $K$ is a compact, connected subset of
    the complex plane. 

 Recall that $K$ is
   called \emph{locally connected} if every point $z\in K$ 
   has arbitrarily small
   connected (relative) neighborhoods in $K$.
   Following Milnor \cite[Chapter 17]{jackdynamicsthird}, 
    we say that $K$ is locally connected \emph{at a point
    $z_0\in K$} if 
    $z_0$ has arbitrarily small connected neighborhoods in $K$.
    Sometimes this property is instead referred to as 
      ``connected im kleinen'', and the term ``locally connected at $z_0$''
      is instead reserved for what Milnor calls
      ``openly locally connected'', see below. 

   Observe
     that local connectivity of $K$ at $z_0$ is equivalent to the condition
     that  the connected component of
       $\{z\in K: |z-z_0|\leq \delta\}$
    containing $z_0$ is a neighborhood of $z_0$ in $K$ for all $\delta>0$.

 \subsection*{Remarks on the definition of local connectivity at a point}
  We note that there are many equivalent definitions of local connectivity 
   of the entire space $K$ which result in different concepts when considered
   only near a point $z_0$. For example, we might say that 
   $K$ is \emph{openly locally connected} at $z_0$ if the point has
   arbitrarily small \emph{open} neighborhoods in $K$. It is well-known
   that $K$ is locally
   connected if and only if it is openly locally connected at every point.
   However, if $K$ is not locally connected as a whole, then it is quite 
   possible that there are points  $z_0\in K$ where the set is 
   locally connected but not openly
   locally connected. (A famous example is the
   ``witch's broom''; see Figure \ref{fig:witchesbroom}.)
 
  Similarly, local connectivity of $K$ implies that $K$ is locally 
   arc-connected, but it is clearly 
   possible for a compact space $K$ to be locally
   connected at a point $z_0$, but not to contain \emph{any} nontrivial 
   curves passing through $z_0$. (See Figure \ref{fig:pathlc}.)

  Finally, 
   Carath\'eodory's theorem is often phrased as a statement on local
   connectivity of the 
   \emph{boundary} $\partial U=\partial K$, rather than of $K$.
   Indeed, these two are equivalent when we consider the entire space; however,
   local connectivity of $\partial U$ at
   $z_0$
   is a strictly stronger condition than local connectivity of $K$ at $z_0$. 
   (See Figure \ref{fig:boundary}.)

  We believe that, in our context, 
   local connectivity of $K$ at $z_0$
   is the most
   natural among the possible concepts to consider. 
   This is vindicated by the 
   fact that we are able to obtain
   natural characterizations of this notion. Also, we should point out
   that our choice places the fewest restrictions on the point $z_0$,
   so that Theorem \ref{thm:main} takes its strongest form with this
   definition.

 \subsection*{Separation theorems and preliminaries}
  We say that two points $z,w\in\Ch$ are \emph{separated} by
   a set $K$ if they belong to different components of $\Ch\setminus K$.
   Similarly, if $U\subset\Ch$ is a domain, we sometimes say that
   $z$ and $w$ are \emph{separated by $K$ in $U$} if they belong to different
   components of $U\setminus K$. 

  We repeatedly use the following standard separation theorem 
   \cite[p.\ 110]{newmanplanetopology}
   due to 
   Janiszewski:
   if
   $K_1$ and $K_2$ are compact subsets of the sphere whose intersection
   is connected, then a pair of points which is not separated by either of
   $K_1$ and $K_2$ is also not separated by the union $K_1\cup K_2$. 

  We will also invoke the \emph{boundary bumping theorem}
   \cite[Theorem 5.6]{continuumtheory}: if $E$ is a subset of
   a compact, connected metric space $K$, then the boundary of every connected
   component of $E$ intersects the boundary of $E$ (in $K$). 

  Let us furthermore remind the reader that a \emph{crosscut} $C$ of a simply
   connected domain $U\subset\Ch$ 
   is a closed Jordan arc which intersects $\Ch\setminus U$
   exactly in its two endpoints. Every crosscut separates $U$ into
   precisely two components. (Since $U$ is homeomorphic to the complex
   plane, this is
   an immediate consequence of the Jordan Curve Theorem; compare
   \cite[Proposition 2.12]{pommerenke}. Observe that the argument
   applies more generally to 
   any injective curve in $U$ which accumulates at $\partial U$ in
   both directions; we use this fact in the Appendix.) 
  
 Finally, we note the following simple result. 

  \begin{lem}[Curves in a subdomain] \label{lem:separationcurve}
   Let $V\subset U$ be a domain, and let 
    $z_0\in \partial U$. Suppose that 
    $\dist(z_0,\partial V\cap U) > \eps$. 
    
   If $w_1,w_2\in V$ can be connected by a curve 
    $\gamma\subset U$ with $\dist(\gamma,z_0) > \eps$, then such a curve
    also exists in $V$.
  \end{lem}
  \begin{proof}
   Let us set $A := \Ch\setminus V\supset K$ and
    $B := K\cup \cl{\D_{\eps}(z_0)}$. By assumption, neither $A$ nor $B$
    separate $w_1$ and $w_2$. We claim that $X := A\cap B$ is connected.

    Indeed, we have $X = K\cup ((U\setminus V)\cap \cl{\D_{\eps}(z_0)})$.
     Suppose, by contradiction, that there was a component
     $L$ of $X$ other than the one containing $K$; then in particular
     $z_0\notin L$. Pick 
     some boundary point $w$
     of $L$ relative to $\cl{\D_{\eps}(z_0)}$. Then $w\in \partial V$, but 
     because $w\in L\subset U$ and $|w-z_0|\leq\eps$, this contradicts
     our assumptions. 

  So $X$ is connected. By Janiszewski's theorem, 
     $A\cup B = (\C\setminus V) \cup \cl{\D_{\eps}(z_0)}$ 
     does not separate $w_1$ and 
     $w_2$, as desired. 
  \end{proof}
 
 \subsection*{A characterization of local connectivity}

\begin{figure}%
  \subfigure[The single comb\label{fig:singlecomb}]{%
                     \resizebox{.425\textwidth}{!}{%
                               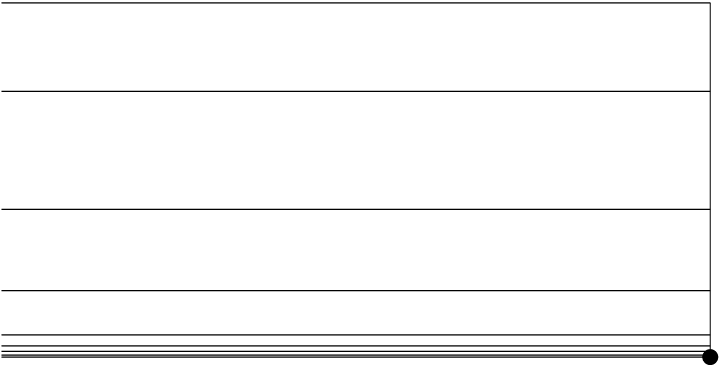}}%
  \hfill
  \subfigure[A modified comb\label{fig:singlecomb_mod}]{%
                        \resizebox{.425\textwidth}{!}{%
                               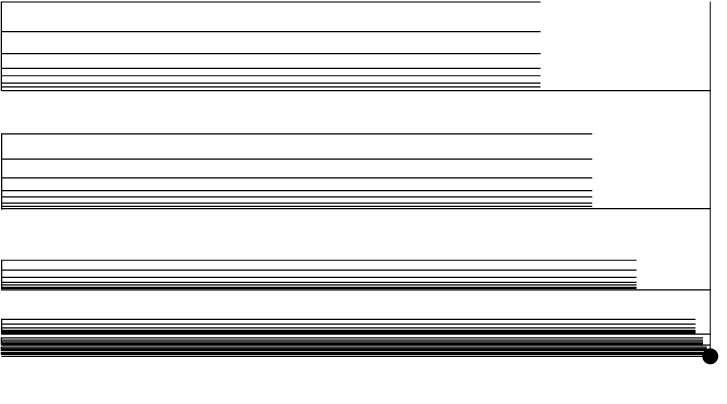}}\\%
  \subfigure[Proposition \ref{prop:localconnectivity}\label{fig:criterion_illustration}]%
                     {\resizebox{.425\textwidth}{!}{%
                               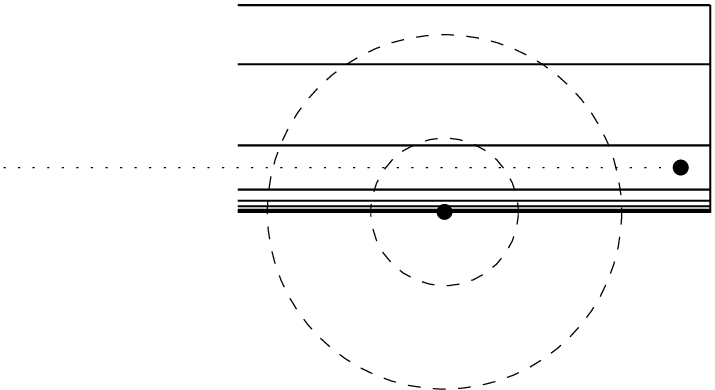}}%
  \hfill
  \subfigure[Proof of Theorem \ref{thm:ursellyoung}\label{fig:ursellyoung}]{%
                      \resizebox{.425\textwidth}{!}{%
                               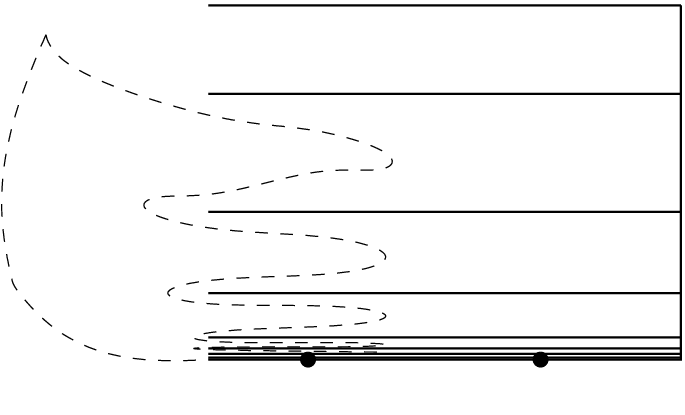}} \caption{%
  The single comb and variations.
    In (a), the set $K=\Ch\setminus U$
    is locally connected at $z_0$ (but not in any other point of
    the interval at the bottom edge of the picture). In contrast,
    the point $z_0$ is not a point of local connectivity in
    (b): this shows that the structure of a prime end's impression,
    and the order of priority from Definition \ref{defn:priority}, 
    is not sufficient to detect local connectivity. Figure (c)
    illustrates the statement of Proposition
    \ref{prop:localconnectivity}: $K$ is not locally connected at
    $z_0$, and any curve connecting $w$ to infinity must intersect
    the disk $\D_{\eps}(z_0)$. 
    Finally, (d) indicates the setup in the
    proof of Theorem \ref{thm:ursellyoung}: the curve $\Gamma$
    accumulates on $z$ but not on $w$, while $\Gamma_0$ accumulates
    on the set $\Pi(p)$ of principal points. Together with their accumulation
    sets, they separate $w$ from the region $V$.\label{fig:singlecombs}}
\end{figure}

  Our proof of Theorem \ref{thm:characterization} (and, by extension,
   of Theorem \ref{thm:main}) relies on the following necessary and sufficient
   condition for local connectivity of $K$ at $z_0$. Compare
   Figure \ref{fig:criterion_illustration}. 

   \begin{prop}[Characterization of local connectivity]
       \label{prop:localconnectivity}
    Let $z_0\in \partial U$. Then $K$ is locally connected at
     $z_0$ if and only if the following holds: for every $\delta>0$,
     there is $\eps>0$ such that every point 
     $w\in U\setminus \D_{\delta}(z_0)$
     can be connected to $\infty$ by a curve
     $\gamma\subset U\setminus \cl{\D_{\eps}(z_0)}$.
   \end{prop}
   \begin{proof}
    Suppose that $K$ is locally connected at $z_0$, and let $\delta>0$. 
     Let $L$ be the connected component of $K\cap \cl{\D_{\delta}(z_0)}$ 
     containing $z_0$. Then $L$ is a compact, connected
     neighborhood of $z_0$ in $K$. I.e., there is 
     $\eps>0$ such that 
        \[ \cl{\D_{\eps}(z_0)}\cap K \subset L. \]
     Let $w_1,w_2\in U\setminus \D_{\delta}(z_0)$.
      Applying Janiszewski's theorem to $K$ and $L\cup \cl{\D_{\eps}(z_0)}$,
      we see that $K\cup \cl{\D_{\eps}(z_0)}$ does not separate $w_1$ from
      $w_2$, as claimed. 

 \smallskip

     For the converse direction, 
      suppose that $K$ is not locally connected at
      $z_0$. Then there is $\delta>0$ such that
      the connected component $L$ of $A := 
      K\cap \cl{\D_{\delta}(z_0)}$ with $z_0\in L$ is not a neighborhood
      of $z_0$ in $K$.

     Let $\eps>0$. Then there is a point $z \in A\setminus L$ with
      $|z-z_0|\leq \eps$. 
      Write $A=A_0\cup A_1$, where
      $A_0$ and $A_1$ are disjoint compact sets with
      $L\subset A_0$ and $z\in A_1$. By the boundary bumping theorem,
      both $L$ and the component of $A$ containing $z$ intersect 
      $\partial \D_{\delta}(z_0)$; in particular, they both
      intersect the circle $\partial \D_{\eps}(z_0)$. 

     So we can pick an arc $C$ of
      $\partial \D_{\eps}(z_0)\setminus A$ which has
      one endpoint in $A_0$ and the other in
      $A_1$.
      Then $\cl{C}$ is a crosscut of $U$; let $V$ be the component of
      $U\setminus C$ which does not contain $\infty$. We claim that
      $V$ contains some point $w$ with $|z_0-w|\geq \delta$. 
      Indeed, applying Janiszewski's theorem first to the 
      sets $A_0$ and $\cl{C}$ and then to
      $A_0\cup \cl{C}$ and $A_1$, we see that we can connect 
      $\infty$ to any point $\tilde{w}\in V$ 
      by a curve not intersecting $A\cup C$. If $|\tilde{w}-z_0|> \delta$,
      we set $w:=\tilde{w}$. Otherwise
      let $w$ be the last intersection point of this
      curve with $\partial \D_{\delta}(z_0)$; then $w\in V$.
     
     By definition of $V$, any curve $\gamma$
      connecting $\infty$ to $w$ must intersect
      $C$, and hence have $\dist(\gamma,z_0)\leq \eps$, as required.
   \end{proof}

 \section{Prime ends} \label{sec:primeends}
  We refer the reader to 
   \cite[Chapter 17]{jackdynamicsthird}
   for an excellent short treatment of the theory of prime ends, and to 
   \cite[Chapter 2]{pommerenke} for further results. Here, we will
   only
   introduce the basic definitions, and give no proofs. As before, $U$
   is a simply connected domain whose complement $K$ omits more than one
   point, and $\infty\in U$ for simplicity. 

  If $C$ is a crosscut of $U$ with $\infty\notin C$, then $U\setminus C$ has
   exactly one component which does not contain $\infty$; let us denote this
   component by $U_C$. A \emph{prime end} of $U$ is represented by
   a sequence of pairwise disjoint crosscuts $(C_n)$ which satisfy
   $\diam C_n\to 0$ and  $C_{n+1}\subset \cl{U_{C_n}}$. Two such sequences
   $(C_n)$ and $(\tilde{C}_n)$ represent the same prime end if 
   $\tilde{C}_j\subset U_{C_n}$ for all $n$ and all sufficiently large $j$,
   and vice versa. 

  The \emph{impression} of $p$ is defined as
      \[ I(p) := \bigcap_{n} \cl{U_{C_n}} \subset \partial U; \]
  an impression is \emph{trivial} if it consists of a single point.
  The set of \emph{principal points}, $\Pi(p)\subset I(p)$, consists of
  those points which are accumulated on by some sequence of crosscuts 
  representing $p$. 

 There is a natural way to define a topology on
   \[ \breve{U} := U \cup \{p:p\text{ is a prime end of $U$}\} \]
  such that a sequence $w_j\in U$ converges to $p$ if and only if
  $w_j\in U_{C_n}$ for all $n$ and all sufficiently large $j$. With this
  topology, $\breve{U}$ is homeomorphic to the closed unit disk
  $\cl{\D}$. In fact, if $\phi:\D\to U$ is a conformal isomorphism, then
  $\phi$ extends continuously to a homeomorphism $\phi:\cl{\D}\to\breve{U}$. 

 We say that a curve $\Gamma:[0,\infty)\to U$ \emph{converges to $p$} if
  $\Gamma(t)\to p$ in the topology of $\breve{U}$. Note that the set of
  accumulation points of such a curve $\Gamma$ necessarily contains $\Pi(p)$.

 \begin{proof}[Proof of Theorem \ref{thm:characterization}]
  Suppose that $K$ is locally connected at $z_0$, and consider a prime end
   $p$ with $z_0\in I(p)$. Let
   $(w_j)$ be a sequence as in the definition of strong minimality
   (Definition \ref{defn:strongminimality}), and
   set $\delta := \inf |z_0-w_j|$. 
   By Proposition   \ref{prop:localconnectivity}, there is $\eps_0>0$
   such that each 
   $w_j$ can be connected to $\infty$ by a curve $\gamma_j$ with 
   $\dist(\gamma_j,z_0)>\eps_0$. 

  Let $(C_n)$ be a sequence of crosscuts representing $p$, and set
   $U_n := U_{C_n}$. Without loss of generality, we may assume that
   no $C_n$ has $z_0$ as an endpoint (otherwise, we simply remove this
   crosscut from the sequence). 

  For sufficiently large $j$, we have
   $w_j\in U_n$, and hence
   $\gamma_j\cap C_n\neq\emptyset$. Since $\dist(\gamma_j,z_0)>\eps_0$ and
   $\diam(C_n)\to 0$, we see that
       \[ \eps_1 := \inf_n \dist(z_0,C_n) > 0; \]
   we set $\eps := \min(\eps_0,\eps_1)$. 

  By construction, we can connect
   $w_j$ and $w_{j+1}$ by a curve $\Gamma_j$ which does not intersect
   $\cl{\D_{\eps}(z_0)}$ (e.g., $\Gamma_j = \gamma_j\cup \gamma_{j+1}$). 
   If $w_j$ and $w_{j+1}$ both belong to the same
   $U_n$, we can, by Lemma \ref{lem:separationcurve}, furthermore choose
   $\Gamma_j$ such that
   $\Gamma_j\subset U_n$. The curve $\Gamma := \bigcup \Gamma_j$ 
   converges to the prime end $p$, contains all points $w_j$, and does not
   intersect $\cl{\D_{\eps}(z_0)}$. 

 \smallskip

  Now suppose, conversely, that $K$ is not locally connected at $z_0$.
   Then, by Proposition \ref{prop:localconnectivity}, there is a constant
   $\delta>0$ with the following property: for every $n\in\N$
   there is a point $\omega_n\in U$ with $|\omega_n-z_0|\geq\delta$ such that
   any curve connecting $\omega_n$ to $\infty$ within $U$ must pass within
   distance at most $1/n$ of $z_0$. 

  Clearly we have $\omega_n\to \partial U$. Let $p$ be some accumulation
   point of $\omega_n$ in the prime end compactification $\breve{U}$ of
   $U$. We define a sequence $(w_j)_{j\geq 0}$ 
   by
   setting $w_0:=\infty$ and
   $w_j := \omega_{n_j}$, where $(\omega_{n_j})_{j\geq 1}$ is a subsequence
   converging to $p$ in $\breve{U}$.

  Then any curve $\Gamma$ tending to $p$ and containing all points
   $w_j$ will, in particular, connect $w_0=\infty$ to $w_j = \omega_{n_j}$.
   Hence $\dist(\Gamma,z_0)=0$ by choice of $(\omega_n)$, which implies
   that $z_0\in I(p)$ and that $z_0$ is not strongly minimal in $p$. 
 \end{proof}

 \begin{cor}[Local connectivity at principal points] \label{cor:principal}
   Suppose that $K$ is locally connected at a principal point $z_0$ of
    $p$. Then $I(p)=\{z_0\}$. 
 \end{cor}
 \begin{proof}
  We prove the contrapositive. Suppose that $I(p)\neq\{z_0\}$; then
   we can pick some $z\in I(p)\setminus \{z_0\}$. 
   Take a sequence of points $w_k\in U$ converging to $z$ such
   that $w_j$ converges to $p$ in the prime end topology. Since $z_0$ is
   a principal point, \emph{any} curve converging to $p$ (containing all $w_j$
   or not) must accumulate on
   $z_0$, and hence $z_0$ is not strongly minimal. By Theorem
   \ref{thm:characterization}, $K$ is not locally connected at $z_0$.
 \end{proof} 
 \begin{remark}
  This fact is easy to prove also without Theorem \ref{thm:characterization}: 
   Let $(C_n)$ be 
   a sequence of crosscuts converging to $z_0$. Pick a closed
   connected neighborhood $Z_{\eps}$ of $z_0$ in $K$ satisfying
   $Z_{\eps}\subset \D_{\eps}(z_0)$; we may suppose without loss of
   generality that $\Ch\setminus Z_{\eps}$ is connected.
   Then, for sufficiently large $n$,
   the curve $C_n$ is a crosscut of $\Ch\setminus Z_{\eps}$.
   It follows that $U_n$, and hence $I(p)$, is contained in
   $\cl{\D_{\eps}(z_0)}$. 
   Since $\eps$ was arbitrary, the claim follows.
 \end{remark}

 \begin{cor}[Carath\'eodory's theorem] \label{cor:caratheodory}
  The set
   $K$ is locally connected if and only if every prime end impression is
   trivial. 
 \end{cor} 
 \begin{proof}
  Since $K$ is locally connected at every point of its interior, 
   we only need to consider local connectivity at points of $\partial U$. 

   If $K$ is \emph{not} locally connected at some point $z_0\in\partial U$, 
    then by Theorem \ref{thm:characterization}
    there is some prime end $p$ with $z_0\in I(p)$ such that $z_0$ is not
    strongly minimal in $p$. In particular, there is some sequence
    $(w_j)$ converging to $p$ but not converging to
    $z_0$, so $I(p)$ is not trivial. 

  On the other hand, if there is some nontrivial prime end impression $I(p)$,
   then $K$ is not locally connected at any point of $\Pi(p)$
   by the previous corollary. 
 \end{proof}

 \section{The left and right wings of a prime end} \label{sec:ursellyoung}
 The article \cite{ursellyoung}
   studied the
   \emph{left and right wings} 
    of a prime end. Today these
    are more commonly referred to as the 
    \emph{left and right cluster sets}, or as the
    \emph{one-sided impressions} \cite{pommerenkecarmonacontinua}. 
    We prefer to use
    ``wing'' here, as it seems to be the original term used
    when these 
    sets are investigated as topological objects related 
    to the domain $U$ and its boundary, rather than 
    as an aspect of a conformal mapping $\phi:\D\to U$.  
    For simplicity, 
    we will nonetheless not
    give a purely topological definition, but rather
    use the conformal mapping $\phi$. 

   Let $p$ be a prime end of $U$.
    We say that a curve $\Gamma:[0,\infty)\to U$ converges to $p$
    \emph{from the left} if $\Gamma$ converges to $p$, and furthermore 
      \[ \im\left(\frac{\phi^{-1}(\Gamma(t))}{p}\right)\geq 0 \]
    for all sufficiently large $t$.
    (Here we again identify the prime end $p$ with the corresponding
     point $p=e^{2\pi i \theta}$ on the unit circle.)
     We say that any accumulation point $z_0\in I(p)$
     of such a curve $\Gamma$ belongs to the \emph{left wing} of
     $p$, and write $I^+(p)$ for all such points. 

     The \emph{right wing} $I^-(p)$ is defined analogously.
      We note that $I(p) = I^+(p)\cup I^-(p)$ and
      $\Pi(p)\subset I^-(p)\cap I^+(p)$. 

 \begin{example}
  In Figure \ref{fig:doublecomb}, the interval at the top of the picture
   is the only nontrivial prime end impression. The midpoint $m$ of this
   interval is the unique principal point; the left and right
   wings are the intervals $[z_0,m]$ and $[m,z_1]$, respectively. 

  In Figures \ref{fig:singlecomb} and \ref{fig:singlecomb_mod}, the 
   prime end at the bottom
   of the picture has a trivial right wing, 
   containing only the unique principal point,
   while the left wing consists of the entire interval at the bottom of
   the picture.
 \end{example}

  The prime end $p$ is called
      \emph{symmetric} if $I^-(p)=I^+(p)$.
    By the Collingwood Symmetry Theorem \cite[Proposition 2.21]{pommerenke}, 
    all but countably many prime ends
    are symmetric. Compare \cite{pommerenkecarmonasymmetric} for   
    interesting results on  symmetric 
    prime ends  (among other things).

  \begin{defn}[Priority] \label{defn:priority}
   Let $z,w\in I^{-}(p)$. 
    We say that 
    \emph{$z$ has priority over $w$} 
    (in $I^-(p)$) if every curve $\Gamma$ 
    which converges to $p$ from the left and 
    accumulates on $w$ must also accumulate on $z$.

   (Priority in $I^+(p)$ is defined analogously.) 
  \end{defn}

  \begin{example}
   In the left wing of the nontrivial prime end expression of 
    Figure \ref{fig:singlecomb}, 
    the order of priority coincides with horizontal order:
    if $z,w$ belong to the interval at the bottom of the picture and $z$ is
    to the left of $w$, then $z$ has priority over $w$. The same is true
    in Figure \ref{fig:singlecomb_mod}. 
  \end{example}

 \begin{thm}[Priority is a total relation \cite{ursellyoung}] 
   \label{thm:ursellyoung}
   Let $z,w$ belong to the same wing of the prime end $p$. Then
    $z$ has priority over $w$ or $w$ has priority over $z$.
 \end{thm}

 Now suppose that $z_0\in I(p)$
  is strongly minimal in $p$. Then $z_0$ is minimal with 
  respect to priority; indeed, $z_0$ cannot have
  priority over any point in either wing. 
  Hence Theorem \ref{thm:ursellyoung} implies:

 \begin{cor}[At most two strongly minimal points] \label{cor:stronglyminimal}
   Each wing of the prime end $p$ contains at most one point which
    is strongly minimal in $p$. \qedoutsideproof
 \end{cor}

 \begin{proof}[Proof of Theorem \ref{thm:main}]
  We just proved that $I(p)$ contains at most two points 
   which are strongly minimal (at most one for each wing). 
   We also know by Theorem
   \ref{thm:characterization} 
   that local connectivity of $K$ at $z_0$ requires
   strong minimality of $z_0$ in $p$. This proves the theorem.

  If the prime end is symmetric, then both wings are equal and
   Corollary \ref{cor:stronglyminimal} implies that $I(p)$ contains at most
   one strongly minimal point, establishing the remark after the statement
   of Theorem \ref{thm:main}.
 \end{proof}

\begin{remark}
 One might ask whether strong minimality can be expressed solely
  in terms of the order of priority on $I(p)$; this is not the case.
  Indeed, the first two examples in Figure \ref{fig:singlecombs} both have
  a prime end whose impression is the bottom interval of the picture, 
  and as discussed above the corresponding orders of priority
  coincide. However, in the first figure
  $z_0$ is strongly minimal, while in the second figure it is not. 
\end{remark}

\section*{Appendix: Proof of the Ursell-Young Theorem}

 \begin{proof}[Proof of Theorem \ref{thm:ursellyoung}]
  Suppose that $z$ and $w$ both
   belong to (say) the left wing of the prime end $p$,
   and that $w$ does not have priority over $z$. That is, there
   is a curve $\Gamma : [0,\infty)\to U$ converging to $p$ from the left
   whose accumulation set $A$ contains $z$ but not $w$. For simplicity,
   let us assume that $\Gamma$ is injective. (It is not hard to see that we
   can always find an injective curve with the same accumulation set.
   Alternatively, with minor modifications the proof will also
   apply in the general case.) 
   We need to show that
   $z$ has priority over $w$. We may assume that $z\notin \Pi(p)$, 
   as every principal point has priority over all other points by definition.

  Let $\Gamma_0:[0,\infty)\to U$ be the ``central curve''
    \[ \Gamma_0(t) := \phi(e^{-1/t}\cdot p) \]
   separating the left and right wings of $p$. The limit set of $\Gamma_0$
   is the set of \emph{radial limit points} of $\phi$ at $p$, which is
   well-known to consist exactly of $\Pi(p)$ \cite[Theorem 2.16]{pommerenke}. 
   We may assume without
   loss of generality that $\Gamma(0)=\Gamma_0(0)$ and that
   $\Gamma \bigl((0,\infty)\bigr)\cap \Gamma_0 = \emptyset$
   (recall that $\Gamma$ converges to $p$ from the left).

  Since $\Gamma_0$ and $\Gamma_1$ accumulate only on $\partial U$,
   the set $U\setminus (\Gamma_0\cup \Gamma)$ has exactly two components. 
   Since $\Gamma_0$ and $\Gamma$ both converge to the same 
   prime end $p$, one of these components, call it $V$, accumulates on
   no other prime ends in the topology of $\breve{U}$. 
            (Compare Figure \ref{fig:ursellyoung}.) We pick
   some arbitrary base point $x\in V$.
  \begin{claim}[Claim 1]
   We have $\partial V = \Gamma_0\cup \Gamma \cup A=: F$. In particular,
   $F$ separates $x$ from $w$.
  \end{claim}
  \begin{subproof}
   Let $\tilde{U}$ be the component of $\Ch\setminus A$ containing $U$ and
    let $\tilde{V}$ be the component of 
    $\tilde{U}\setminus (\Gamma_0\cup \Gamma)$ containing $V$. Then
    $V=\tilde{V}\cap U$; we need
    to show that $V=\tilde{V}$. 

   Indeed, otherwise $\tilde{V}$ is a neighborhood of some
    point of $\partial U$; 
    in particular $\tilde{V}$ and hence $V$
    contains some crosscut of $U$. However,
    this is a contradiction to the choice of $V$, since every crosscut
    accumulates on two distinct prime ends of $U$. 
  \end{subproof} 

   \begin{claim}[Claim 2]
     Let $\eps$ be sufficiently small. Let $T>0$ be minimal with
      $|\Gamma(T) - z|\leq \eps$. Then the set
       \[ K_1(\eps) := \cl{\D_{\eps}(z)}\cup \Gamma\bigl([T,\infty)\bigr)\cup
                        A \]
      does not separate $x$ and $w$ (in $\Ch$). 
   \end{claim}
   \begin{subproof}
     Let $w'\in U$ with $|w - w'| < \dist(w, A)$. By connecting 
      $x$ to $w'$ in $U$ and $w'$ to $w$ by a straight line segment,
      we obtain a curve $\alpha\subset \C\setminus A$ which connects
      $x$ to $w$. Set 
      $\delta := \dist(\alpha,A)$. 
      If $\eps$ is sufficiently small, then $K_1(\eps)$ is contained in 
      a $\delta$-neighborhood of $A$, and hence does not intersect $\alpha$.
   \end{subproof}

  Let $(C_n)$ be a sequence of crosscuts representing $p$;
   we may choose these so that $C_n\cap \Gamma_0$ consists of a single
   point for every $n$. Let $U_n$ be the component of
   $U\setminus (C_n\cup \Gamma_0)$ which contains $\Gamma(t)$ for large $t$.
   Then a curve converges to $p$ from the left
   if and only if it is eventually contained
   in every $U_n$. 

  \begin{claim}[Claim 3]
   Let $\eps>0$. Then there are $n_0$ and $\delta$ with the following
    property.
    If $n\geq n_0$ and $w'\in U_n$ with $|w-w'| < \delta$, then
    any curve in $U_n$ connecting $w'$ to $C_n$ intersects
    $\cl{\D_{\eps}(z)}$.
  \end{claim}
  \begin{subproof}
    By decreasing $\eps$, if necessary, we may assume that $x$
     can be connected to every crosscut $C_n$ by a curve 
    $\beta_n \subset V$ which does not intersect
    $\cl{\D_{\eps}(z)}$. This is possible because
    $\Gamma_0$ does not accumulate on $z$. 

   By Claim 2, we may
    also assume that $K_1 = K_1(\eps)$ does not separate $x$ and $w$. 
    Consider the set
     \[ K_2 := \cl{\D_{\eps}(z)} \cup A \cup \Gamma_0 \cup
                 \Gamma\bigl([0,T]\bigr). \]
    Observe that $K_1\cap K_2 = \D_{\eps}(z) \cup A$
      and
      $F\subset K_1\cup K_2$. Hence it follows from Claim 1
      and Janiszewski's
      theorem that $K_2$ separates $x$ and $w$.

    Choose $\delta$ sufficiently small that 
    $\D_{\delta}(w)\cap K_2 = \emptyset$, 
    and choose $n_0$ such that $C_n\cap \cl{\D_{\eps}(z)}=\emptyset$
    and $U_n\cap \Gamma\bigl([0,T]\bigr)=\emptyset$ for $n\geq n_0$. 

   If 
    $w'\in U_n$ with $n\geq n_0$ and $|w'-w|<\delta$, then 
    $K_2$ separates $w'$ and $x$. 
    Let $\gamma\subset U_n$ be a curve
    connecting $w'$ to $C_n$. Combining $\gamma$ with 
    a piece of $C_n$ and 
    the curve $\beta_n$, we obtain a curve 
    in $U\setminus (\Gamma_0 \cup \Gamma([0,T]))$ 
    connecting $w'$ to
    $x$. This curve must intersect $K_2$, and hence
    $\cl{\D_{\eps}(z)}$. Since $C_n$ and $\beta_n$ do not intersect
    this disk, it follows that $\gamma$ intersects
    $\cl{\D_{\eps}(z)}$, 
    as claimed.
  \end{subproof}

  The proof of Theorem \ref{thm:ursellyoung} is now complete: suppose
   that $\tilde{\Gamma}$ is
   a curve converging to $p$ from the left and accumulating on $w$, and
   let $\eps>0$. Let $\delta$ and $n_0$ be as in Claim
   3, and pick $n\geq n_0$ sufficiently large so that
   $\tilde{\Gamma}\not\subset U_n$. Since $\tilde{\Gamma}$ contains
   some point $w'\in U_n$ with $|w'-w| < \delta$, it follows that
   $\tilde{\Gamma}$ intersects $\cl{\D_{\eps}(z)}$.
 \end{proof}

\bibliographystyle{hamsplain}
\bibliography{../../Biblio/biblio}

\end{document}